\documentclass[letter, 12pt]{amsart}
\usepackage[utf8]{inputenc}
\usepackage{xcolor}
\usepackage{hyperref}
\hypersetup{colorlinks=true,linkcolor={blue}}
\usepackage[inner=2.3cm,outer=2.3cm, bottom=3.3cm]{geometry}

\usepackage{amsmath, mathtools}
\usepackage{cleveref}
\usepackage{amssymb}
\usepackage{amsfonts}
\usepackage{amsthm, enumitem}

\newtheorem{thm}{\bf Theorem}[section]

\theoremstyle{definition}
\newtheorem{definition}[thm]{\bf Definition}
\newtheorem{strategy}[thm]{\bf Strategy}

\theoremstyle{remark}

\numberwithin{equation}{section}

\DeclareMathOperator{\gr}{gr}
\DeclareMathOperator{\hgt}{ht}

\def \fm{\mathfrak{m}}
\def \fn{\mathfrak{n}}
\def \fp{\mathfrak{p}}
\def \G{\mathcal{G}}

\title[Computing multiplicity sequences]{Computing multiplicity sequences} 

\author[Justin Chen]{Justin Chen}
\address{School of Mathematics, Georgia Institute of Technology, Atlanta, GA 30332}
\email{jchen646@gatech.edu}

\author[Youngsu Kim]{Youngsu Kim}
\address{Department of Mathematics, California State University San Bernardino, San Bernardino, CA 92407}
\email{youngsu.kim@csusb.edu}

\author[Jonathan Monta\~no]{Jonathan Monta\~no$^{*}$}
\address{School of Mathematical and Statistical Sciences, Arizona State University, P.O. Box 871804, Tempe, AZ 85287-18041}
\email{montano@asu.edu}
\thanks{$^{*}$ The third  author is  supported by NSF Grant DMS \#2001645/2303605.}

\begin{document}
	
	\begin{abstract}
		The {\tt MultiplicitySequence} package for \texttt{Macaulay2} computes the multiplicity sequence of a graded ideal in a standard graded ring over a field, as well as several invariants of monomial ideals related to integral dependence. 
		We discuss two strategies implemented for computing multiplicity sequences: one via the bivariate Hilbert polynomial, and the other via the technique of general elements.
	\end{abstract}
	
	\maketitle
	
	\section{Introduction}
	Let $(R,\fm,k)$ be a $d$-dimensional Noetherian local ring with maximal ideal $\fm$ and residue field $k$. Let $I$ be an $R$-ideal. 
	If $I$ is $\fm$-primary, then the \emph{Hilbert-Samuel multiplicity} of $I$ is defined as the degree of the standard graded algebra 
	$$\gr(I)=\bigoplus_{n=0}^\infty I^n/I^{n+1}$$
	i.e.\ the normalized leading coefficient of its Hilbert polynomial.
	This classical numerical invariant has been the base of several important results in commutative algebra and algebraic geometry.
	For example, a classical result of Rees states that the Hilbert-Samuel multiplicity gives an effective criterion for deciding whether two ideals have the same integral closure, provided $R$ is formally equidimensional \cite{Rees61}.
	Rees' Theorem is of fundamendal importance in singularity theory as it is a key component in the proof of Teissier's \emph{Principle of Specialization of Integral Dependence} (PSID), which provides a fiberwise numerical criterion for a family of hypersurfaces with isolated singularities to be equisingular \cite{Tes}.
	
	The \emph{$j$-multiplicity} and \emph{$\varepsilon$-multiplicity} are extensions of the Hilbert-Samuel multiplicity to arbitrary ideals.
	These multiplicities were originally introduced in \cite{AM} and \cite{KV}, respectively, in large part to extend Rees' Theorem to the non-$\fm$-primary case.
	Such extensions were obtained in \cite{FM} and \cite{UV}, but with the requirement of having to localize at all prime ideals of $R$.
	
	The {\it multiplicity sequence} of an arbitrary ideal $I$ in $R$ is a sequence of $d+1$ non-negative integers corresponding to the leading coefficients of the second sum transform of the bivariate Hilbert polynomial of the standard bigraded algebra 
	\begin{equation} \label{grgr}
		\G := \gr(\fm\gr(I))=\bigoplus_{i,j=0}^{\infty}\frac{\fm^iI^j+I^{j+1}}{\fm^{i+1}I^j+I^{j+1}}.
	\end{equation} 
	This sequence is a particular case of the multiplicities defined by Kleiman and Thorup in \cite[\S 8]{KT}, and it was also considered  by Gaffney and Gassler \cite{GG} in the analytic case, and by Achilles and Manaresi \cite{AM2} in our general setting. 
	In the recent work \cite{PTUV}, Polini, Trung, Ulrich, and Validashti extended Rees' Theorem to arbitrary ideals without the need for localizations, by using multiplicity sequences. 
	More precisely, they show that if $R$ is formally equidimensional, then ideals $I\subseteq J$ have the same integral closure if and only if their multiplicity sequences agree (the forward direction was previously obtained by Ciuperc\u{a} in \cite{Ciu}).
	Furthermore, the authors of \cite{PTUV} develop a PSID for arbitrary ideals using the multiplicity sequence, demonstrating the importance of this invariant. 
	
	The main goal of the {\tt MultiplicitySequence} package in \texttt{Macaulay2} \cite{GS} is to compute the multiplicity sequence of graded ideals in standard graded rings over a field.
	Two strategies have currently been implemented for doing so: the first one is based on the definition via the bivariate Hilbert polynomial of $\G$, and the second strategy is based on the technique of general elements, cf. \cite[Theorem 4.1]{AM2}. 
	Finally, the package also includes a number of methods related to multiplicities and integral dependence, which have been adapted to the case of monomial ideals. 
	
	\section{Multiplicity Sequence}
	Throughout, we keep the same notation as in the introduction. 
	For a module $M$, $\lambda(M)$ denotes the length of $M$.
	
	\subsection{Associated bi-graded ring}
	The second sum transform of the \emph{bivariate Hilbert polynomial} of $\G$ is the polynomial $P(m,n)$ that agrees with
	\begin{equation}\label{lengthGij}
		h(m,n)=\sum_{i=0}^{m}\sum_{j=0}^{n}\lambda(\G_{i,j}),\quad \text{where}\quad \G_{i,j}=\frac{\fm^iI^j+I^{j+1}}{\fm^{i+1}I^j+I^{j+1}}
	\end{equation}
	
	for $m,n\gg 0$.
	The polynomial $P(m,n)$ can be written in the form
	$$P(m,n)=\sum_{i=0}^d\frac{c_i(I)}{(d-i)!i!}m^{d-i}n^i+\text{(lower degree terms)}$$
	with $c_i(I) \in \mathbb{Z}_{\ge 0}$ for $i=0,\ldots, d$ \cite{VDW}.  
	
	\begin{definition}
		The sequence $c_0(I), \ldots, c_d(I)$ is called the \emph{multiplicity sequence} of $I$.
	\end{definition}
	One has $c_i(I)=0$ if $i<d-\dim R/I$ or $i>\ell(I)$, where $\ell(I):=\dim \gr(I)\otimes_R k$ is the {\it analytic spread} of $I$ \cite[Proposition 2.3]{AM2}.
	Moreover, $c_{d}(I)$ equals the $j$-multiplicity of $I$. 
	In particular, if $I$ is $\fm$-primary, then $c_d(I)$ is the Hilbert-Samuel multiplicity of $I$ while $c_i(I) = 0$ for $i \neq d$.
	
	For purposes of \texttt{Macaulay2} computation, we take the local ring $R$ to be of the form $A_{\fn}$, where $A$ is a standard graded algebra over a field and $\fn$ is its irrelevant ideal (note that lengths of graded modules do not change under localizing at $\fn$).
	We now describe our first strategy for computing the multiplicity sequence.
	
	\begin{strategy}\label{str1}
		Given an ideal $I$, we compute the bigraded algebra $\G$ using 
		{\tt tangentNormalCone} (which iteratively calls {\tt normalCone}).
		Subsequently, the method {\tt hilbertSequence} extracts the relevant coefficients of the Hilbert polynomial $P(m,n)$ of $\G$ from the Hilbert series of $\G$.
	\end{strategy}
	
	\Cref{str1} is the default strategy for computing the multiplicity sequence, and is executed whenever \texttt{multiplicitySequence} is called without specifying any options.
	We illustrate its use in the following example:
	
	\begin{verbatim}
		Macaulay2, version 1.17
		i1 : needsPackage "MultiplicitySequence";
		i2 : S = QQ[a..e]/(ideal(a-b,c)*ideal(c,d,e));
		i3 : I = ideal"a2-bd,b4,e3";
		i4 : multiplicitySequence I
		o4 = HashTable{2 => 3 }
		               3 => 12
		i5 : hilbertSequence tangentNormalCone I
		o5 =     0 1 2 3   
		       +----------
		     3 | . . 3 12 
		     2 | . . 2 .  
		     1 | . 1 . .  
		     0 | . . . .
	\end{verbatim}
	
	In the output \texttt{o4} above, the multiplicity sequence is displayed as a hash table, indicating that $c_2(I) = 3$ and $c_3(I) = 12$.
	The coefficients of the Hilbert polynomial of $\G$ are displayed in \texttt{o5} as a 2-dimensional table, whose top row is precisely the multiplicity sequence of $I$.
	
	The most time-consuming step in \Cref{str1} is that of computing (a presentation of) $\G$ -- the Hilbert series and coefficient extraction are comparatively fast.
	For convenience, this expensive step is cached upon completion, so later calls to \texttt{multiplicitySequence} for a given ideal are nearly instant.
	
	\subsection{General Elements}
	Our second strategy is based on \Cref{thmGenEl} below which uses the method of general elements.
	For a local ring $S$, we denote by $e(S)$ the Hilbert-Samuel multiplicity of its maximal ideal. 
	
	\begin{thm}[{\cite[Remark 2.3]{PTUV}}]\label{thmGenEl}
		Suppose $R$ is equidimensional and catenary with infinite residue field.
		For any $i\geqslant 0$ and general elements $x_1,\ldots, x_i$ of $I$, one has
		\begin{equation}\label{multiSeqGen}
			c_i(I)=
			\sum_{\fp}\lambda\left(\frac{R_\fp}{(x_1,\ldots, x_{i-1})R_\fp:I^\infty +x_iR_\fp}\right)e(R/\fp),
		\end{equation}
		where the sum ranges over the set of prime ideals
		\begin{equation}\label{setOfPrimes}
			\{\fp\in V(I)\mid \hgt \fp = i,\, \fp\supset (x_1,\ldots, x_{i-1}):I^\infty \}
		\end{equation}
		and by convention the colon ideal $ (x_1,\ldots, x_{i-1}):I^\infty$ is 0 if $i=0$ and is $0:I^\infty$ if $i=1$.
	\end{thm}
	
	In view of \Cref{thmGenEl}, one could compute $c_i(I)$ by choosing general elements $x_1, \ldots, x_\ell \in I$, and then computing the various lengths and multiplicities in \Cref{multiSeqGen}.
	However, this necessitates localizing at all the primes $\mathfrak{p}$ appearing above, which is undesirable for \texttt{Macaulay2} computation. Thus we take a different approach, as explained below.
	
	\begin{strategy}\label{str2}
		Via \Cref{thmGenEl}, we identify $c_i(I)$ with $e(R/J_i)$ for a suitable $R$-ideal $J_i$, and the latter can be computed in \texttt{Macaulay2} using a combination of {\tt degree} and {\tt normalCone} (in particular, avoiding localizations).
		The ideal $J_i$ is constructed as follows: first, compute the minimal primes of the ideal $(x_1,\ldots, x_{i-1}):I^\infty + (x_i)$. 
		Next, set $K$ to be the intersection of these minimal primes that do not contain $I$. 
		Finally, define $J_i := ((x_1,\ldots, x_{i-1}):I^\infty + x_i) \colon K^\infty$.
		To see that $c_i(I) = e(R/J_i)$, note that we may identify the set of primes (\ref{setOfPrimes}) as 
		\begin{align*}
			\{\fp\in V(I) &\mid \hgt \fp = i,\, \fp\supset (x_1,\ldots, x_{i-1}):I^\infty \} \\
			&= \{\fp\in V((x_1,\ldots, x_{i-1}):I^\infty + x_i) \mid \hgt \fp = i \} \cap V(I) \\
			&= \{\fp\in V((x_1,\ldots, x_{i-1}):I^\infty + x_i) \mid \hgt \fp = i \} \setminus (\operatorname{Spec} (R) \setminus V(I) )\\
			&= \{\fp\in V(J_i) \mid \hgt \fp = i \}.
		\end{align*}
		Then by the associativity formula for Hilbert-Samuel multiplicity, the sum in \Cref{multiSeqGen}, taken over the last set above, is precisely $e(R/J_i)$. 
		
	\end{strategy}
	
	We illustrate the use of \Cref{str2}.
	Note that the index $i$ for $c_i(I)$ is specified here.
	\begin{verbatim}
		i6 : R = QQ[a..d];
		i7 : I = ideal "a2,ab,b3,ad - bc,c2-bd";
		i8 : multiplicitySequence(3, I, Strategy => "generalElements") -- c_3(I)
		o8 = 5
		i9 : multiplicitySequence(4, I, Strategy => "generalElements") -- c_4(I)
		o9 = 7
	\end{verbatim}
	
	For specific values of $i$, \Cref{str2} may be faster than \Cref{str1}.
	However, for computing the entire multiplicity sequence, \Cref{str1} tends to outperform \Cref{str2}, hence our choice of \Cref{str1} as the default strategy.
	
	As noted before, for $i=d$, the coefficient $c_d(I)$ is equal to the $j$-multiplicity of $I$, which has been studied by several authors, see e.g.\ \cite{NU,PX,JM,JMV}. We isolate this case in the method {\tt jMultiplicity}, which is based on code written by H.\ Schenck and J.\ Validashti.
	
	\section{Methods for Monomial Ideals}
	
	Monomial ideals carry combinatorial structure which often allows for special algorithms.
	The {\tt MultiplicitySequence} package contains a few methods dedicated to monomial ideals, such as \texttt{newtonPolyhedron, monomialReduction}, and specialized algorithms for \texttt{analyticSpread} and \texttt{jMultiplicity} in the case of monomial ideals. 
	These methods utilize the Newton polyhedron of a monomial ideal and scale much more efficiently than general methods.
	For comparison, we show the difference in timings for some of these methods:
	
	\begin{verbatim}
		i10 : I = monomialIdeal"ab2,bc3,cd4,da5";
		i11 : elapsedTime jMultiplicity I^3
		-- 0.874315 seconds elapsed
		o11 = 9639
		i12 : elapsedTime jMultiplicity ideal I^3
		-- 456.039 seconds elapsed
		o12 = 9639
		i13 : elapsedTime analyticSpread I^5
		-- 0.515529 seconds elapsed
		o13 = 4
		i14 : elapsedTime analyticSpread ideal I^5
		-- 42.4524 seconds elapsed
		o14 = 4
	\end{verbatim}
	
	\section*{Acknowledgments}
	We thank the anonymous referees for their careful reading of the paper and Macaulay2 package. The second author thanks D.\ Eisenbud, D.\ Grayson, and M.\ Stillman for organizing a {\tt Macaulay2} day during the special year in commutative algebra 2012-2013 at MSRI, where he learned how to write a package.  The third  author is  supported by NSF Grant DMS \#2001645/2303605.

\end{document}